\documentclass[12pt,thmsa]{article}
\usepackage{amsmath}
\usepackage{graphicx}
\usepackage{amsfonts}
\usepackage{amssymb}
\usepackage{mathrsfs}

\begin{document}

\begin{center}
{\Large \textbf{Robustifying multiple-set linear canonical analysis with  S-estimator}}

\bigskip

Ulrich  DJEMBY BIVIGOU,  Guy Martial  NKIET

\bigskip

Universit\'{e} des Sciences et Techniques de Masuku,  Franceville, Gabon.

\bigskip

E-mail adresses : dastakreas@yahoo.fr,  guymartial.nkiet@mathsinfo.univ-masuku.com.

\bigskip
\end{center}

\noindent\textbf{Abstract.}We consider  a robust version of multiple-set linear canonical analysis obtained by using a S-estimator of  covariance operator. The related influence functions are derived. Asymptotic properties of  this robust method are obtained and  a robust test for mutual non-correlation is introduced. 

\bigskip



\section{Introduction}

\noindent
Letting \ ($\Omega, \mathcal A, P)$ \  be a probability space$,$ and \ $K$ \ be an integer such that \ $K \geq 2,$ \ we consider random variables \ $X_1 ,  \cdots  , X_K$ with values in Euclidean vector spaces \ $\mathcal{X}_1 ,  \cdots  , \mathcal{X}_K$ \ respectively. We suppose that for any $k\in\{1,\cdots,K\}$, one has $\mathbb{E}(X_k)=0$ and  $\mathbb E(\left\|X_k\right\|_k^2)$ $<\infty ,$ where \ $\left\| \cdot \right\|_k$ \ denotes the norm induced by the inner product   $\left\langle \cdot\  ,\  \cdot\right\rangle_{k}$  of $\mathcal{X}_k$.  The multiple-set  linear canonical analysis (MSLCA)  of $X_1,\cdots,X_K$   is a statistical method that permits to analyse the relationships among these variables. It has been introduced for many years (e.g., \cite{gifi}) and has been extensively studied (e.g., \cite{gardner},  \cite{hwang}, \cite{nkiet},  \cite{takane},\cite{tenenhaus}). Formally, considering the  random variable $X=(X_1,\cdots, X_K)$ with values into the space  $\mathcal{X}:=\mathcal{X}_1\times \mathcal{X}_2\times\cdots\mathcal{X}_K$,  MSLCA  is the search of    a sequence \ $\left(\alpha^{(j)}\right)_{ 1\leq j \leq q}$ \  of vectors of $\mathcal{X}$, where $q =$ dim$(\mathcal{X})$,  satisfying:
\begin{equation}\label{acg}
\alpha^{(j)} =\text{arg}\ \underset{\alpha \in \mathcal C_j}{\text{max}}\  \mathbb E\left(< \ \alpha\  ,\  X\ >^2_{\mathcal{X}}\right) ,
\end{equation}
\noindent
where 
$\mathcal C_1 = \ \left\{  \alpha \in  \mathcal{X} / \  \sum_{k=1}^K \ var\left(< \ \alpha_k\  ,\  X_k\  >_{k}\right) = 1 \right\}$  and,  for $ j \geq 2$:
\[
\mathcal C_j = \ \left\{  \alpha \in  \mathcal C_1 / \  \sum_{k=1}^K \ cov\left(< \ \alpha^{(r)}_k\  ,\  X_k\  >_{k} , < \ \alpha_k\  ,\  X_k\ >_{k}\right) = 0 , \ \forall \ r \in \left\{ 1 , \cdots  , j - 1 \right\}\right\}.
\]
A solution of the above maximization problem is obtained from the  spectral
analysis of a given  operator $T$. For  $(k,\ell) \in \left\{ 1 , \cdots  , K \right\}^2 $,   considering  the covariance operators
$V_{k\ell} = \mathbb E(X_\ell\otimes X_k) = V_{\ell k}^*$   and \ $V_{k} := V_{kk}$, 
where \ $\otimes$ \  denotes the tensor product such that \ $x\otimes y$\  is the linear map :  $h \mapsto < x, h > y$,   letting  $ \tau_k$  be the canonical projection
$\tau_k: \alpha =(\alpha_1,\cdots,\alpha_K)\in \mathcal{X} \mapsto \alpha_k \in \mathcal{X}_k$ and assuming that each $V_k$ is invertible, we have $T= \Phi^{-1/2}\Psi\Phi^{-1/2}$, where 
$
\Phi =  \sum_{k=1}^K \tau_k^*V_{k}\tau_k$ and $\Psi =   \sum_{k=1}^{K} \sum_{\underset{\ell\neq k}{\ell=1}}^{K} \tau_k^*V_{k\ell}\tau_\ell$.
If $\left\{ \beta^{(1)}, \cdots, \beta^{(q)}\right\}$	 is an orthonormal basis of $\mathcal{X}$ \  such that $\beta^{(j)}$  is an eigenvector of $T$ associated with the $j$-th largest eigenvalue $\rho_j,$ then we obtain a solution of (1) by taking
$\alpha^{(j)}= \Phi^{-1/2}\beta^{(j)}$. Classical estimation of MSLCA is based on   empirical covariance operators that are known to be  very sensitive to outliers. This makes this method a non robust one  and highlights the interest of providing a robust estimation of MSLCA as it was done for others multivariate statistical methods such as principal component analysis, discriminant analysis, canonical correlation analysis (\cite{crouxdehon}, \cite{crouxfilm}, \cite{crouxhaes}, \cite{taskinen}).  A natural way for doing that is to replace the covariance operator of $X$ by a robust estimator. Among such robust estimators, the S-estimator has been extensively studied (\cite{davies}, \cite{lopuhaa}, \cite{lopuhaa2}) and it is known to have good robustness and efficiency properties.  In this paper,  we propose a robust version of MSLCA based on S-estimator of the covariance operator.  This estimation procedure is introduced in Section 2. The related influence functions are derived in Section 3. Section 4 is devoted to the asymptotic properties of the introduced estimator, and a robust test for mutual non-correlation is proposed in Section 5.

\section{Estimation of MSLCA based on S-estimator}
\noindent  We assume that the following condition holds:
\bigskip

\noindent $(\mathscr{A}_1):$   $X$ has  an elliptical contoured distribution with density
$
f_X(x) = (\textrm{det}(V))^{-1/2}h(<x , V^{-1}x >_{\mathcal{X}})
$, where $h : [0, +\infty[ \rightarrow [0, +\infty[$ is a function having a strictly negative derivative $h'$. 

\bigskip

\noindent Let $\{X^{(1)},\cdots,X^{(n)}\}$  be an i.i.d. sample of   $X$, we consider a fixed real   $b_0$ and a function $\xi : \mathbb{R} \rightarrow \mathbb{R}$. We denote by $\mathcal{P}\left(\mathcal{X}\right)$ the set  of positive definite symmetric operators from $\mathcal{X}$ to itself. The  S-estimators $\widetilde{\mu}_{n}$   and   $\widetilde{V}_{n}$ of the mean and the covariance operator of $X$    respectively are given by the pair $(\widetilde{\mu}_{n},\widetilde{V}_{n})$ that  minimizes  the determinant  det$\left(G\right)$ over all $(\mu,G) \in\mathcal{X}\times\mathcal{P}\left(\mathcal{X}\right)$   that satisfy 
$
 \frac{1}{n} \sum_{i=1}^{n} \xi\left(\left\|G^{-1/2}(X^{(i)} - \mu)\right\|_{\mathcal{X}}\right)\leq b_0 . 
$
\noindent
It is well known that these estimators are robusts and have high breakdown points (see, e.g., \cite{davies}). From them,  we can introduce an estimator of MSLCA which is expected to be also robust. Indeed,  putting $
\widetilde{\Phi}_n =  \sum_{k=1}^{K} \tau^*_k\widetilde{V}_{k.n}\tau_k$ and $\widetilde{\Psi}_n =   \sum_{k=1}^{K} \sum_{\underset{\ell\neq k}{\ell=1}}^{K}  \tau^*_k\widetilde{V}_{k\ell.n}\tau_\ell$, 
where $\widetilde{V}_{k.n} = \tau_k\widetilde{V}_{n}\tau^*_k$  and  $\widetilde{V}_{k\ell.n} = \tau_k\widetilde{V}_{n}\tau^*_\ell$,  we estimate  $T$   by
$
\widetilde{T}_n =  \widetilde{\Phi}^{-1/2}_n \widetilde{\Psi}_n\widetilde{\Phi}^{-1/2}_n 
$. 
Considering the eigenvalues  $\widetilde{\rho}_{1.n} \geq \widetilde{\rho}_{2.n} \geq \cdots\geq \widetilde{\rho}_{q.n}$ of  $\widetilde{T}_n$   and  $\left\{ \widetilde{\beta}^{(1)}_n \cdots  \widetilde{\beta}^{(q)}_n \right\}$ an orthonormal basis of $\mathcal{X}$  such that \ $\widetilde{\beta}^{(j)}_n$ \  is an eigenvector of  $\widetilde{T}_n$    associated with $\widetilde{\rho}_{j.n},$ we estimate $\rho_{j}$ \  by  $\widetilde{\rho}_{j.n},$  $\beta^{(j)}$ \  by \ $\widetilde{\beta}^{(j)}_n$ \  and \ $\alpha^{(j)}$ \  by \ $\widetilde{\alpha}^{(j)}_n = \widetilde{\Phi}^{-1/2}_n  \widetilde{\beta}^{(j)}_n$.  

\section{Influence functions}
\noindent For studying the effect of a small amount of  contamination at a given point    on  MSLCA it is important,  as usual in robustness litterature (see \cite{hampel}), to use influence function. More precisely, we have to derive expressions of the  influence functions related to the functionals that give $T$, $\rho_j$ and $\alpha^{(j)}$ (for $1\leq j \leq q$) at  the distribution $\mathbb{P}_X$ of $X$. Recall that the influence function of a functional $S$ at $\mathbb{P}$ is defined as
\[
\textrm{IF}\left(x;S,\mathbb{P}\right)=\lim_{\varepsilon\downarrow 0}\frac{S\left((1-\varepsilon)\mathbb{P}+\varepsilon\delta_x\right)-S(\mathbb{P})}{\varepsilon},
\]
where $\delta_x$ is the Dirac measure putting all its mass  in $x$.
In order to derive the influence functions related to the above estimator of MSLCA$,$ we have to specify the functional that corresponds to it. We impose the following properties on the loss function  $\xi$:

\bigskip

\noindent $(\mathscr{A}_2)$:  $\xi$  is symmetric,   has a continuous derivative  $\psi $  and is such that  $\xi(0)=0$;

\bigskip

\noindent $(\mathscr{A}_3)$:  there exists  $c_0 > 0$ such that \ $\xi$  is strictly increasing on \ $[0 ,\  c_0]$ \ and constant  on    $[c_0 ,\  +\infty[$;

\bigskip 

\noindent $(\mathscr{A}_4)$:  the function $t \mapsto \psi(t)t^{-1}$ is continuous and  bounded.

\bigskip

\noindent  For example, the   function  $\xi (t)=\frac{c^2}{6}\left(1-\left(1-\frac{t^2}{c^2}\right)^3\right)\,\textrm{\textbf{1}}_{[-c,c]}(t)+\frac{c^2}{6}\,\textrm{\textbf{1}}_{\mathbb{R}\backslash [-c,c]}(t)$, where $c>0$, satisfies the above conditions. Its derivative is the Tukey's biweight function   $\psi (t)=\left(1-\frac{t^2}{c^2}\right)^2\,\textrm{\textbf{1}}_{[-c,c]}(t)$.
The functional $\mathbb V_s$  related to the aforementioned S-estimator of V is defined in \cite{lopuhaa2}  (see also \cite{lopuhaa}); it is such that   $\mathbb V_s(\mathbb{P})$ is  the solution to the problem of minimizing the determinant  det$\left(G\right)$ over all $\mu \in \mathcal{X}$ and $ G \in \mathcal{P}( \mathcal{X})$  that satisfy 
\begin{eqnarray*}
 \int_{\mathcal{X}} \xi\left(\left\|G^{-1/2}(x - \mu)\right\|_{\mathcal{X}}\right)\,d\mathbb{P}(x)&\leq& b_0 . 
\end{eqnarray*}
\noindent
It is known that at elliptical distribution  $\mathbb V_s(\mathbb P_X) = V$ (see \cite{lopuhaa},  p.222). Therefore,  the functional $\mathbb T_s$   defined as $\mathbb T_s(\mathbb P) = f\left(\mathbb V_s(\mathbb P)\right)^{-1/2} g\left(\mathbb V_s(\mathbb P)\right) f\left(\mathbb V_s(\mathbb P)\right)^{-1/2}$,  
where $f\left(A\right) = \displaystyle\sum_{k=1}^K \tau_k^*\tau_k A\tau_k^*\tau_k$  and  $ g\left(A\right) = \displaystyle \sum_{k=1}^{K} \sum_{\underset{\ell\neq k}{\ell=1}}^{K} \tau_k^*\tau_kA\tau_\ell^*\tau_\ell $, we have 	$\mathbb T_s(\mathbb{P}_X)=T$ .  Putting $T_s = \mathbb T_s(\mathbb P_X) $,
\begin{equation}\label{l}
\lambda(x,V)=\sum_{k=1}^K \sum_{\underset{\ell\neq k}{\ell=1}}^K  -\frac{1}{2}\tau_k^*\left(x_k\otimes x_k\right) V_{k\ell}\tau_\ell - \frac{1}{2} \tau_\ell^* V_{\ell k}\left(x_k\otimes x_k\right)\tau_k+\tau_k^* \left(x_\ell\otimes x_k\right)\tau_\ell
\end{equation}
\[
\gamma_1 = \frac{2\pi^{q/2}}{\Gamma(q/2)(q + 2)} \int_{0}^{+\infty}\left(    \psi'(r)r^2 + (q + 1)\psi(r)r\right) r^{q - 1}h(r^2)dr   
\]

and 
\[
\gamma_2 =  \frac{2\pi^{q/2}}{\Gamma(q/2)} \int_{0}^{+\infty}\psi(r)r^{q} h(r^2)dr , 
\]
   $\Gamma$ \ being the usual gamma function,  we have:
\bigskip

\noindent
\textbf{Theorem 3.1.}\textsl{ We suppose that the assumptions $(\mathscr{A}_1)$ to $(\mathscr{A}_4)$ hold. Then}
\begin{eqnarray*}
\textrm{IF}(x; T_s, \mathbb{P}_X) &=&   \frac{q}{\gamma_1} \psi\left(\left\|V^{-1/2}x\right\|_{\mathcal{X}}\right)\,\,\left\|V^{-1/2}x\right\|_{\mathcal{X}}^{-1}\,\lambda(x,V).
\end{eqnarray*}

\noindent From the properties of $\psi$, it is easily seen that $\textrm{IF}(x; T_s, \mathbb{P}_X) $ equals  $0$ if  $\Vert V^{-1/2}x\Vert_\mathcal{X}>c_0$. Otherwise, we have $\Vert x\Vert_\mathcal{X}\leq c_0\,\Vert V^{1/2}\Vert_\infty$, where 
    $\left\|\cdot\right\|_{\infty}$ denotes  the usual operators norm defined by $\left\|A\right\|_{\infty}=\text{sup}_{x\neq 0 }\left(\left\|Ax\right\| /\left\|x\right\| \right)$.  Then, it is easy to check the inequality
\[
\underset{x \in \mathcal{X}}{\sup}\left\|\textrm{IF}(x; T_s,\mathbb{P}_X)\right\|_{\infty} \leq  \sup_{t \in \mathbb R^*_{+}}\left(  \frac{\psi(t)}{t}\right)  \frac{Kq}{\left|\gamma_1\right|}\left(K - 1\right)\left( \left\|V\right\|_{\infty} + 1 \right)c_0^2\left\|V^{1/2}\right\|^2_{\infty}
\]
that  shows that the influence function  is bounded and, therefore, that  the estimation procedure is robust. Now, we give the influence functions related to the canonical coefficients and the canonical directions obtained from the robust MSLCA introduced above. For $j\in\{1,\cdots,q\}$, denoting by $\mathbb{R}_{s\cdot j}$ (resp. $\mathbb{B}_{s\cdot j}$; resp. $\mathbb{A}_{s\cdot j}$ ) the functional such that $\mathbb{R}_{s\cdot j}(\mathbb{P})$ is the $j$-th largest eigenvalue of $\mathbb{T}_{s}(\mathbb{P})$ (resp. the associated eigenvector; resp. $\mathbb{A}_{s\cdot j}(\mathbb{P})= f\left(\mathbb{V}_{s}(\mathbb{P} )\right)^{-1/2}\mathbb{B}_{s\cdot j}(\mathbb{P})$ ), we put  $\rho_{s\cdot j}=\mathbb{R}_{s\cdot j}(\mathbb{P}_X)$,  $\beta^{(j)}_{s}=\mathbb{B}_j(\mathbb{P}_X)$  and $\alpha^{(j)}_{s}=\mathbb{A}_{s\cdot j}(\mathbb{P}_X)$.  
Putting
\begin{eqnarray}\label{h}
H\left(\xi , \psi , V , x\right) &= &  \frac{2}{\gamma_2}\left(\xi\left(\left\|V^{-1/2}x\right\|_{\mathcal{X}}\right) - b_0\right)\mathbb{I} \nonumber\\
& &+   \frac{q}{\gamma_1}\psi\left(\left\|V^{-1/2}x\right\|_{\mathcal{X}}\right)\left\|V^{-1/2}x\right\|_{\mathcal{X}}\left(  \frac{1}{\left\|V^{-1/2}x\right\|^2_{\mathcal{X}}} -   \frac{1}{q}\right)\mathbb{I}
\end{eqnarray} 
and
\begin{eqnarray}\label{lj}
\lambda_j( x , V) 
&=&\sum_{k=1}^{K} \sum_{\underset{l\neq k}{\ell=1}}^{K}\sum_{\underset{m\neq j}{m=1}}^{q} \frac{1}{\rho_j - \rho_m}\bigg(<\beta_k^{(m)} ,  x_k  >_k< x_\ell , \beta_\ell^{(j)}  >_\ell \nonumber\\
& &-\frac{1}{2}  < \beta_k^{(m)} ,  x_k >_k < x_k , V_{k\ell} \beta_\ell^{(j)}  >_k \nonumber\\  
& &- \frac{1}{2} < x_k ,  V_{k\ell }\beta_\ell^{(m)}  >_\ell < x_k , \beta_k^{(j)} >_k \bigg)\beta^{(m)}\\
& &-\frac{1}{2}\left( \sum_{k=1}^{K}\left[ \tau_k^*\left(x_k\otimes x_k\right)\tau_k + < \beta_k^{(j)}, x_k>_k^2\mathbb{I}\right]  -  2\mathbb{I}\right)\beta^{(j)} \nonumber,
\end{eqnarray}
where $\mathbb{I}$ denotes the identity operator of $\mathcal{X}$, we have:
\bigskip

\noindent
\textbf{Theorem 3.2.} \textsl{We suppose that the assumptions the assumptions $(\mathscr{A}_1)$ to $(\mathscr{A}_4)$ hold.  Then,  for
any  $j \in  \left\{1, \cdots , q\right\},$  we have:}
\\\\
$\left(i\right)$ \ $\textrm{IF}(x; \rho_{s.j}, \mathbb P_X) =   \frac{q}{\gamma_1}\psi\left(\left\|V^{-1/2}x\right\|_{\mathcal{X}}\right)\left\|V^{-1/2}x\right\|^{-1}_{\mathcal{X}}  \sum_{k=1}^{K} \sum_{\underset{l\neq k}{\ell=1}}^{K} < \beta_k^{(j)}  , x_k >_{k} < x_\ell - V_{\ell k}x_k  ,  \beta_\ell^{(j)}  >_{\ell}$.

\bigskip

\noindent
$\left(ii\right)$ \textsl{If $\rho_1>\cdots>\rho_q$, then } $\textrm{IF}( x ; \alpha_s^{(j)} ,\mathbb P_X) =   \frac{q}{\gamma_1}\psi\left(\left\|V^{-1/2}x\right\|_{\mathcal{X}}\right)\left\|V^{-1/2}x\right\|^{-1}_{\mathcal{X}}\,\lambda_j( x , V)  - H\left(\xi , \psi , V , x\right)\beta^{(j)}  $.

\section{Asymptotic distributions}
\noindent
We first establish asymptotic normality for $\widetilde{T}_{n}$. Putting 
\[
\beta_3 =  \frac{2\pi^{q/2}}{\Gamma(q/2)} \int_{0}^{+\infty} \frac{4}{q + 2}\psi(r)r^{q + 2} h'(r^2)dr ,
\]

\noindent
we have:

\bigskip

\noindent
\textbf{Theorem 4.1. }\textsl{We suppose that the assumptions  $(\mathscr{A}_1)$ to $(\mathscr{A}_4)$ hold  and that $\mathbb{E}\left(\Vert X\Vert_\mathcal{X}^4\right)<+\infty$. Then,   $ \sqrt{n}\left(\widetilde{T}_n - T\right)$ converges in distribution,  as 
$n \rightarrow +\infty$,  to a random variable $U_s$ having a normal distribution in $\mathcal L(\mathcal{X}),$ with mean $0$  and 
covariance operator equal to that of the random operator}
\begin{eqnarray}\label{zs}
Z_s &=&   -2q\beta^{-1}_3\psi\left(\left\|V^{-1/2}X\right\|_{\mathcal{X}}\right)\,\left\|V^{-1/2}X\right\|_{\mathcal{X}}^{-1}\,\, \lambda(X,V).
\end{eqnarray}
\noindent
This result allows to consider a robust  test for mutual non-correlation, that is the test  for the hypothesis  $\mathscr{H}_0 : \forall (k, \ell) \in \left\{1, ... , K \right\}^2,$  $k \neq \ell,$  $V_{k\ell} = 0$
against the alternative  
$\mathscr{H}_1 : \exists (k, \ell) \in  \left\{1, ... , K\right\}^2$,   $k \neq \ell$,   $V_{k\ell} \neq 0$. We take as test statistic the random variable $
\widetilde{S}_{n} = \displaystyle \sum_{k=2}^{K} \sum_{\ell=1}^{k - 1}\textrm{tr}\left( \pi_{k\ell}\left(\widetilde{T}_{n}\right)\pi_{k\ell}\left(\widetilde{T}_{n}\right)^*\right)$,   where $\pi_{k\ell}$ is the operator $A\mapsto \tau_k A\tau_\ell^\ast$.
Then, putting 
\begin{eqnarray}
\kappa_0 &=& \frac{-2\beta^{-1}_3}{(q + 1)}\mathbb E\bigg(\psi\left(\left\|X\right\|_{\mathcal X}\right)\left\|X\right\|_{\mathcal X}^3\bigg),
\end{eqnarray}
\noindent
we have:
\bigskip

\noindent
\textbf{Theorem 4.2. }\textsl{We suppose that the assumptions   $(\mathscr{A}_1)$ to $(\mathscr{A}_4)$ hold  and that $\mathbb{E}\left(\Vert X\Vert_\mathcal{X}^4\right)<+\infty$. Then,  under
$\mathscr{H}_0$,  the sequence $(\kappa_0)^{-1}n\widetilde{S}_{n}$ converges in distribution$,$ as $n \rightarrow +\infty,$ to \ $\chi_d^2,$ \ where $d = \sum_{k=1}^{K}\sum_{l=1}^{k-1} p_kp_l$  with \ $p_k = \textrm{dim}(\mathcal X_k).$}

\bigskip

\noindent In practice,   $\kappa_0$ is replaced  by a consistent estimator; for example by $\widehat{\kappa}_0= \frac{-2\beta^{-1}_3}{n(q + 1)}\sum_{i=1}^n \psi\left(\left\|X^{(i)}\right\|_{\mathcal X}\right)\left\|X^{(i)}\right\|_{\mathcal X}^3.$
 \section{Proof of Theorem 3.1}

\noindent
It is shown in  Lopuha$\ddot{\textrm{a}}$ (1989)  (see  Corollary 5.2,  p. 1672) that under spherical distribution $\mathbb P^0_X,$ one has
\begin{eqnarray*}
\textrm{IF}(x; V_s, \mathbb P^0_X) &=&   \frac{2}{\gamma_2}\left(\xi(\left\|x\right\|_{\mathcal{X}}) - b_0\right)\mathbb I +   \frac{q}{\gamma_1}\psi(\left\|x\right\|_{\mathcal{X}})\left\|x\right\|_{\mathcal{X}}\left(  \frac{x\otimes x}{\left\|x\right\|^2_{\mathcal{X}}} -   \frac{1}{q}\mathbb I\right) .
\end{eqnarray*} 
\noindent
Then, affine equivariant property implies that under  elliptical model given in assumtion $(\mathscr{A}_1)$ we have:
\begin{eqnarray}\label{ifvs}
\textrm{IF}(x; V_s, \mathbb P_X) &=& V^{1/2}\left\{  \frac{2}{\gamma_2}\left(\xi\left(\left\|V^{-1/2}x\right\|_{\mathcal{X}}\right) - b_0\right)\mathbb I\right. \nonumber\\
& &+\left.   \frac{q}{\gamma_1}\psi(\left\|V^{-1/2}x\right\|_{\mathcal{X}})\left\|V^{-1/2}x\right\|_{\mathcal{X}}\left(  \frac{(V^{-1/2}x)\otimes (V^{-1/2}x)}{\left\|V^{-1/2}x\right\|^2_{\mathcal{X}}} -   \frac{1}{q}\mathbb I\right)\right\}V^{1/2} \nonumber\\
&=&   \frac{2}{\gamma_2}\left(\xi\left(\left\|V^{-1/2}x\right\|_{\mathcal{X}}\right) - b_0\right)V \nonumber\\
& &+   \frac{q}{\gamma_1}\psi\left(\left\|V^{-1/2}x\right\|_{\mathcal{X}}\right)\left\|V^{-1/2}x\right\|_{\mathcal{X}}\left(  \frac{x\otimes x}{\left\|V^{-1/2}x\right\|^2_{\mathcal{X}}} -   \frac{1}{q}V\right).
\end{eqnarray} 
\noindent
Putting  $V_s = \mathbb V_{s} \left(\mathbb P_X\right) = V,$  we have  $f\left(\mathbb V_{s} \left(\mathbb P_X\right)\right) = f\left(V\right) = \mathbb I.$
Thus

\begin{eqnarray*}
\mathbb T_s \left(\mathbb P_{\varepsilon,x}\right) - \mathbb T_s \left(\mathbb P_X\right) &=& f\left(\mathbb V_{s} \left(\mathbb P_{\varepsilon,x}\right)\right)^{-1/2}g\left(\mathbb V_{s} \left(\mathbb P_{\varepsilon,x}\right)\right)f\left(\mathbb V_{s} \left(\mathbb P_{\varepsilon,x}\right)\right)^{-1/2} - g\left(V_s\right)\\
&=& \left(f\left(\mathbb V_{S} \left(\mathbb P_{\varepsilon,x}\right)\right)^{-1/2} - \mathbb I\right)g\left(\mathbb V_{s} \left(\mathbb P_{\varepsilon,x}\right)\right)f\left(\mathbb V_{s} \left(\mathbb P_{\varepsilon,x}\right)\right)^{-1/2}\\
&&+\left(g(\mathbb V_{s} \left(\mathbb P_{\varepsilon,x}\right) - V_s)\right)f\left(\mathbb V_{s} \left(\mathbb P_{\varepsilon,x}\right)\right)^{-1/2}\\
&&+g\left(V_s\right)\left(f\left(\mathbb V_{s} \left(\mathbb P_{\varepsilon,x}\right)\right)^{-1/2} - \mathbb I\right).
\end{eqnarray*} 
\noindent
where $\mathbb P_{\varepsilon,x} = (1 - \varepsilon)\mathbb P_X + \varepsilon \delta_x$  with  $\varepsilon \in [0; 1]$. Then$,$ using the equality
\begin{eqnarray}\label{deco}
A^{-1/2} - \mathbb I &=& -A^{-1}(A - \mathbb I) \left(A^{-1/2} + \mathbb{I}\right)^{-1} 
\end{eqnarray} 
\noindent
we obtain:
\begin{eqnarray*}
\mathbb T_s \left(\mathbb P_{\varepsilon,x}\right) - \mathbb T_s \left(\mathbb P_X\right) &=& -f\left(\mathbb V_{s} \left(\mathbb P_{\varepsilon,x}\right)\right)^{-1/2}f\left(\mathbb V_{s} \left(\mathbb P_{\varepsilon,x}\right) - \mathbb V_{s} \left(\mathbb P_X\right)\right)\left(f\left(\mathbb V_{s} \left(\mathbb P_{\varepsilon,x}\right)\right)^{-1/2} + \mathbb I\right)^{-1}\\
&&\times g\left(\mathbb V_{s} \left(\mathbb P_{\varepsilon,x}\right)\right)f\left(\mathbb V_{s} \left(\mathbb P_{\varepsilon,x}\right)\right)^{-1/2}\\
&&+\left(g(\mathbb V_{s} \left(\mathbb P_{\varepsilon,x}\right) - V_S)\right)f\left(\mathbb V_{s} \left(\mathbb P_{\varepsilon,x}\right)\right)^{-1/2}\\
&&-g\left(V_s\right)f\left(\mathbb V_{s} \left(\mathbb P_{\varepsilon,x}\right)\right)^{-1/2}f\left(\mathbb V_{s} \left(\mathbb P_{\varepsilon,x}\right) - \mathbb V_{s} \left(\mathbb P_X\right)\right)\left(f\left(\mathbb V_{s} \left(\mathbb P_{\varepsilon,x}\right)\right)^{-1/2} + \mathbb I\right)^{-1}.
\end{eqnarray*} 
\noindent
Then$,$ from
\begin{eqnarray*}
\textrm{IF}( x ; V_s ,\mathbb P)&=&  \lim_{\varepsilon\rightarrow 0} \frac{\mathbb V_{s} \left(\mathbb P_{\varepsilon,x}\right) - \mathbb V_{s}\left(\mathbb P_X\right) }{\varepsilon}
\end{eqnarray*}
\noindent
and the continuity of the maps  $A \mapsto A^{-1}$ and  $A \mapsto A^{-1/2} ,$ we
deduce that
\begin{eqnarray}\label{ifts}
\textrm{IF}(x; T_s, \mathbb P_X)&=& -  \frac{1}{2}f\left(\textrm{IF}(x; V_s, \mathbb P_X)\right)g\left(V_{s}\right)\nonumber\\
&&+g\left(\textrm{IF}(x; V_s, \mathbb P_X)\right)\nonumber\\
&&-  \frac{1}{2}g\left(V_{s}\right)f\left(\textrm{IF}(x; V_s, \mathbb P_X)\right).
\end{eqnarray} 
\noindent
Since  $f\left(V\right) = \mathbb I ,$ and under elliptical model  $g\left(V_{s}\right) = g\left(V\right) ,$ inserting (\ref{ifvs}) in (\ref{ifts}) gives the equality
\begin{eqnarray*}
\textrm{IF}(x; T_s, \mathbb P_X)&=&   \frac{q}{\gamma_1}\psi\left(\left\|V^{-1/2}x\right\|_{\mathcal{X}}\right)\left\|V^{-1/2}x\right\|^{-1}_{\mathcal{X}}\bigg\{ -  \frac{1}{2}f\left(x\otimes x\right)g\left(V_s\right)  \\
& &-  \frac{1}{2}g\left(V_s\right)f\left(x\otimes x\right)  + g\left(x\otimes x\right)\bigg\}\\
&&+ \left\{  \frac{2}{\gamma_2}\left(\xi\left(\left\|V^{-1/2}x\right\|_{\mathcal{X}}\right) - b_0\right) -   \frac{q}{\gamma_1}\psi\left(\left\|V^{-1/2}x\right\|_{\mathcal{X}}\right)\left\|V^{-1/2}x\right\|_{\mathcal{X}}\right\}\\
&& \times \left\{ -  \frac{1}{2}f\left(V\right)g\left(V_s\right)  -  \frac{1}{2}g\left(V_s\right)f\left(V\right)  + g\left(V\right)\right\}.
\end{eqnarray*}
Using properties of the tensor product, is is easy to check that
\[
 -  \frac{1}{2}f\left(x\otimes x\right)g\left(V\right)  -  \frac{1}{2}g\left(V\right)f\left(x\otimes x\right)  + g\left(x\otimes x\right)=\lambda(x,V).
\]
Hence
\begin{eqnarray*}
\textrm{IF}(x; T_s, \mathbb P_X)&=&   \frac{q}{\gamma_1}\psi\left(\left\|V^{-1/2}x\right\|_{\mathcal{X}}\right)\left\|V^{-1/2}x\right\|^{-1}_{\mathcal{X}}\,\lambda(x,V)\\
&&+ \bigg\{  \frac{2}{\gamma_2}\left(\xi\left(\left\|V^{-1/2}x\right\|_{\mathcal{X}}\right) - b_0\right) \\
& &-   \frac{q}{\gamma_1}\psi\left(\left\|V^{-1/2}x\right\|_{\mathcal{X}}\right)\left\|V^{-1/2}x\right\|_{\mathcal{X}}\bigg\}\left(g\left(V\right)  - g\left(V_s\right)\right)\\
&=&   \frac{q}{\gamma_1}\psi\left(\left\|V^{-1/2}x\right\|_{\mathcal{X}}\right)\left\|V^{-1/2}x\right\|^{-1}_{\mathcal{X}}\,\lambda(x,V).
\end{eqnarray*}

\section{Proof of Theorem 3.2}

$(i)$. \ From Lemma 3 in Croux and Dehon (2002),  we obtain

\begin{eqnarray*}
\textrm{IF}(x; \rho_{s.j}, \mathbb P_X) &=& < \beta^{(j)}  , \textrm{IF}(x; T_s, \mathbb P_X) \beta^{(j)}>_{\mathcal{X}}\\
&=&   \frac{q}{\gamma_1}\psi\left(\left\|V^{-1/2}x\right\|_{\mathcal{X}}\right)\left\|V^{-1/2}x\right\|^{-1}_{\mathcal{X}} < \beta^{(j)}  , \lambda(x,V) \beta^{(j)} >_{\mathcal{X}}\\
&=&   \frac{q}{\gamma_1}\psi\left(\left\|V^{-1/2}x\right\|_{\mathcal{X}}\right)\left\|V^{-1/2}x\right\|^{-1}_{\mathcal{X}}  \sum_{k=1}^{K} \sum_{\underset{\ell\neq k}{\ell=1}}^{K} < \beta_k^{(j)}  , x_k >_{k} < X_\ell - V_{\ell k}x_k  ,  \beta_l^{(j)}  >_{l}.
\end{eqnarray*} 

\noindent
$(ii)$. \ Since  $f\left(V_s\right) = f\left(V\right) = \mathbb I ,$ we obtain by applying the second part of Lemma 3 in Croux and Dehon (2002):

\begin{eqnarray}
\textrm{IF}(x; \beta_s^{(j)}, \mathbb P_X) &=& \sum_{\underset{m\neq j}{m=1}}^{q}   \frac{1}{\rho_{j} - \rho_{m} }< \beta^{(m)}  , \textrm{IF}(x; T_s, \mathbb P_X) \beta^{(j)}>_{\mathcal{X}}\beta^{(m)} \nonumber\\
&& -   \frac{1}{2}< \beta^{(j)}  , \textrm{IF}(x; f(V_s), \mathbb P_X) \beta^{(j)}>_{\mathcal{X}}\beta^{(j)}.
\end{eqnarray} 

\noindent
Using   (\ref{ifvs}), the equalities $\textrm{IF}(x; f(V_s), \mathbb P_X) = f\left(\textrm{IF}(x; V_s, \mathbb P_X)\right)$,  $f(V)=\mathbb{I}$  and   
Theorem 3.1, we obtain:
\begin{eqnarray*}
\textrm{IF}(x; \beta_s^{(j)}, \mathbb P_X) &=&    \frac{q}{\gamma_1}\psi\left(\left\|V^{-1/2}x\right\|_{\mathcal{X}}\right)\left\|V^{-1/2}x\right\|^{-1}_{\mathcal{X}} \lambda(x,V)\sum_{\underset{m\neq j}{m=1}}^{q}   \frac{1}{\rho_{j} - \rho_{m} }< \beta^{(m)}  , \beta^{(j)}>_{\mathcal{X}}\beta^{(m)} \nonumber\\
&& -   \frac{1}{2} \bigg\{\frac{2}{\gamma_2}\left(\xi\left(\left\|V^{-1/2}x\right\|_{\mathcal{X}}\right) - b_0\right)\\
& &+   \frac{q}{\gamma_1}\psi\left(\left\|V^{-1/2}x\right\|_{\mathcal{X}}\right)\left\|V^{-1/2}x\right\|^{-1}_{\mathcal{X}}< \beta^{(j)} ,(f(x\otimes x) - \mathbb I )\beta^{(j)}>_{\mathcal{X}}\\
&& +   \frac{q}{\gamma_1}\psi\left(\left\|V^{-1/2}x\right\|_{\mathcal{X}}\right)\left\|V^{-1/2}x\right\|_{\mathcal{X}}\left(  \frac{1}{\left\|V^{-1/2}x\right\|^2_{\mathcal{X}}} -   \frac{1}{q}\right)\bigg\}\beta^{(j)}\\
 &=&   \frac{q}{\gamma_1}\psi\left(\left\|V^{1/2}x\right\|_{\mathcal{X}}\right)\left\|V^{-1/2}x\right\|^{-1}_{\mathcal{X}}\eta_j(x,V) \nonumber\\
&& -   \frac{1}{2}\bigg\{  \frac{2}{\gamma_2}\left(\xi\left(\left\|V^{-1/2}x\right\|_{\mathcal{X}}\right) - b_0\right)\\
& & +   \frac{q}{\gamma_1}\psi\left(\left\|V^{-1/2}x\right\|_{\mathcal{X}}\right)\left\|V^{-1/2}x\right\|_{\mathcal{X}}\left(  \frac{1}{\left\|V^{-1/2}x\right\|^2_{\mathcal{X}}} -   \frac{1}{q}\right)\bigg\}\beta^{(j)},
\end{eqnarray*} 
where
\begin{eqnarray*}
\eta_j(x,V)&=&\frac{q}{\gamma_1}\psi\left(\left\|V^{-1/2}x\right\|_{\mathcal{X}}\right)\left\|V^{-1/2}x\right\|^{-1}_{\mathcal{X}} \bigg\{\lambda(x,V)\sum_{\underset{m\neq j}{m=1}}^{q}   \frac{1}{\rho_{j} - \rho_{m} }< \beta^{(m)}  , \beta^{(j)}>_{\mathcal{X}}\beta^{(m)} \\
& &-   \frac{1}{2}< \beta^{(j)} ,(f(x\otimes x) - \mathbb I )\beta^{(j)}>_{\mathcal{X}}\beta^{(j)}\bigg\}.
\end{eqnarray*}
\noindent
On the other hand,  since
\[
\alpha_s^{(j)}(\mathbb P_{X}) =  f\left(\mathbb V_{s}(\mathbb P_X)\right)^{-1/2}\beta_s^{(j)}(\mathbb P_X) = f\left(V\right)^{-1/2}\beta_s^{(j)}(\mathbb P_X) = \beta_s^{(j)}(\mathbb P_X),
\]
it follows
\begin{eqnarray*}
\alpha_s^{(j)}(\mathbb P_{\varepsilon,x}) - \alpha_s^{(j)}(\mathbb P_{X}) &=&  f\left(\mathbb V_{s}(\mathbb P_{\varepsilon,x})\right)^{-1/2}\beta_s^{(j)}(\mathbb P_{\varepsilon,x}) - \beta_s^{(j)}(\mathbb P_X)\\
&=& f\left(\mathbb V_{s}(\mathbb P_{\varepsilon,x})\right)^{-1/2}\left(\beta_s^{(j)}(\mathbb P_{\varepsilon,x}) - \beta_s^{(j)}(\mathbb P_X)\right)\\
&& + \left(f\left(\mathbb V_{S}(\mathbb P_{\varepsilon,x})\right)^{-1/2} - \mathbb I\right)\beta_s^{(j)}(\mathbb P_X)\\
\end{eqnarray*} 
\noindent
Then using (\ref{deco}), we obtain:
\begin{eqnarray*}
\alpha_s^{(j)}(\mathbb P_{\varepsilon,x}) - \alpha_s^{(j)}(\mathbb P_{X}) &=&  f\left(\mathbb V_{S}(\mathbb P_{\varepsilon,x})\right)^{-1/2}\beta_s^{(j)}(\mathbb P_{\varepsilon,x}) - \beta_s^{(j)}(\mathbb P_X)\\
&=& f\left(\mathbb V_{s}(\mathbb P_{\varepsilon,x})\right)^{-1/2}\left(\beta_s^{(j)}(\mathbb P_{\varepsilon,x}) - \beta_s^{(j)}(\mathbb P_X)\right)\\
&& - f\left(\mathbb V_{s}(\mathbb P_{\varepsilon,x})\right)^{-1}\left(f(V_{s}(\mathbb P_{\varepsilon,x}) - V_{s}(\mathbb P_{X}))\right)\\
&&\times\left(f\left(\mathbb V_{S}(\mathbb P_{\varepsilon,x})\right)^{-1/2} + \mathbb I\right)^{-1}\beta_s^{(j)}(\mathbb P_X)\\
\end{eqnarray*} 
\noindent
From the continuity of the maps $A \mapsto A^{-1}, A \mapsto A^{-1/2},$  and the equalities
$ \lim_{\varepsilon\rightarrow 0}f\left(\mathbb V_{S}(\mathbb P_{\varepsilon,x})\right) = f\left(V_{s}(\mathbb P_{X})\right)=f\left(V\right) = \mathbb I ,$ we deduce that
\begin{eqnarray*}
\textrm{IF}( x ; \alpha_s^{(j)} ,\mathbb P)&=&\textrm{IF}(x; \beta_s^{(j)}, \mathbb P_X) -   \frac{1}{2}f\left(\textrm{IF}(x; V_s, \mathbb P_X)\right)\beta_s^{(j)}\\
&=&    \frac{q}{\gamma_1}\psi\left(\left\|V^{-1/2}x\right\|_{\mathcal{X}}\right)\left\|V^{-1/2}x\right\|^{-1}_{\mathcal{X}}\left\{\eta_j(x,V) -   \frac{1}{2}\left(f\left(x\otimes x\right)-\mathbb{I}\right)\beta^{(j)}\right\}\nonumber\\
&& +   \frac{1}{2}\bigg\{  \frac{q}{\gamma_1}\psi\left(\left\|V^{-1/2}x\right\|_{\mathcal{X}}\right)\left\|V^{-1/2}x\right\|^{-1}_{\mathcal{X}}\left(f\left(x\otimes x\right)-\mathbb{I}\right)\beta^{(j)} \\
& & - f\left(\textrm{IF}(x; V_s, \mathbb P_X)\right)\beta_s^{(j)} \bigg\}\nonumber\\
&& -   \frac{1}{2}\bigg\{  \frac{2}{\gamma_2}\left(\xi\left(\left\|V^{-1/2}x\right\|_{\mathcal{X}}\right) - b_0\right) \\
& &+   \frac{q}{\gamma_1}\psi\left(\left\|V^{-1/2}x\right\|_{\mathcal{X}}\right)\left\|V^{-1/2}x\right\|_{\mathcal{X}}\left(  \frac{1}{\left\|V^{-1/2}x\right\|^2_{\mathcal{X}}} -   \frac{1}{q}\right)\bigg\}\beta^{(j)}. 
\end{eqnarray*}
\noindent
It is easy to check that  $\eta_j(x,V) -   \frac{1}{2}\left(f\left(x\otimes x\right)-\mathbb{I}\right)\beta^{(j)}=\lambda_j(x,V)$. Then,  since  $\beta_s^{(j)} = \beta^{(j)}$, we deduce from (\ref{ifvs}) that:
\begin{eqnarray*}
\textrm{IF}( x ; \alpha_s^{(j)} ,\mathbb P)&=&   \frac{q}{\gamma_1}\psi\left(\left\|V^{-1/2}x\right\|_{\mathcal{X}}\right)\left\|V^{-1/2}x\right\|^{-1}_{\mathcal{X}}\lambda_j(x,V) \nonumber\\ 
&&+   \frac{q}{2\gamma_1}\psi\left(\left\|V^{1/2}x\right\|_{\mathcal{X}}\right)\left\|V^{-1/2}x\right\|^{-1}_{\mathcal{X}}\left(f\left(x\otimes x\right)-\mathbb{I}\right)\left(\beta^{(j)} - \beta_s^{(j)}\right)\\
&&  -   \frac{1}{2}\bigg\{  \frac{2}{\gamma_2}\left(\xi\left(\left\|V^{-1/2}x\right\|_{\mathcal{X}}\right) - b_0\right)\\
& & +   \frac{q}{\gamma_1}\psi\left(\left\|V^{-1/2}x\right\|_{\mathcal{X}}\right)\left\|V^{-1/2}x\right\|_{\mathcal{X}}\left(  \frac{1}{\left\|V^{-1/2}x\right\|^2_{\mathcal{X}}} -   \frac{1}{q}\right)\bigg\}\left( \beta^{(j)} + \beta_s^{(j)}\right)\\
&=&   \frac{q}{\gamma_1}\psi\left(\left\|V^{-1/2}x\right\|_{\mathcal{X}}\right)\left\|V^{-1/2}x\right\|^{-1}_{\mathcal{X}}\lambda_j(x,V) \nonumber\\ 
&&  - \bigg\{  \frac{2}{\gamma_2}\left(\xi\left(\left\|V^{-1/2}x\right\|_{\mathcal{X}}\right) - b_0\right) \\
& &+   \frac{q}{\gamma_1}\psi\left(\left\|V^{-1/2}x\right\|_{\mathcal{X}}\right)\left\|V^{-1/2}x\right\|_{\mathcal{X}}\left(  \frac{1}{\left\|V^{-1/2}x\right\|^2_{\mathcal{X}}} -   \frac{1}{q}\right)\bigg\}\beta^{(j)}.
\end{eqnarray*}

\section{Proof of Theorem 4.1}
\subsection{A preliminary lemma}
\noindent
The following lemma gives the asymptotic distribution of the random variable
\begin{eqnarray*}
\widetilde{H}_{n} &=&   \sqrt{n}\left(\widetilde{V}_n - V\right).
\end{eqnarray*}

\noindent
\textbf{Lemma 1.} \textit{We suppose that the assumptions   $(\mathscr{A}_1)$ to $(\mathscr{A}_4)$ hold  and that $\mathbb{E}\left(\Vert X\Vert_\mathcal{X}^4\right)<+\infty$. Then,  $\widetilde{H}_{n}$ converges 
in distribution,  as $n \mapsto +\infty$,  to a random variable  having a normal distribution in $\mathcal L(\mathcal{X})$ with mean $0$ and covariance operator equal to that of }
\begin{eqnarray*}
\mathcal{Z} &=& -2q\beta^{-1}_3\psi\left(\left\|V^{-1/2}X\right\|_{\mathcal{X}}\right)\left\|V^{-1/2}X\right\|^{-1} X\otimes X\\
&& - 2   \left(  \frac{\xi(\left\|V^{-1/2}X\right\|_{\mathcal{X}}) - b_0}{q\beta_1} -    \frac{\psi(\left\|V^{-1/2}X\right\|_{\mathcal{X}})\left\| V^{-1/2}X\right\|_{\mathcal{X}}}{\beta_3}\right) V.
\end{eqnarray*}
\\\\
\noindent
\textit{Proof.} Let $\theta$ and $\phi$ be the functions from $\mathbb{R}_+$ to $\mathbb{R}$ defined as $\theta (t)= 2q\beta^{-1}_3\psi(t) t$ and $\phi(t) = 2q^{-1}\beta^{-1}_1\left(\xi(t) - b_0\right) - 2\beta^{-1}_3\psi(t)t $.       Using affine equivariant property$,$ we deduce from the proof of Corollary 2  in  Lopuha$\ddot{\textrm{a}}$ (1997) (see  p. 235)
that:

\begin{eqnarray*}
\widetilde{H}_{n} &=& V^{1/2}\bigg( -   \frac{1}{  \sqrt{n}}  \sum_{i=1}^{n}\bigg\{  \frac{\theta(\left\|V^{-1/2}X^{(i)}\right\|_{\mathcal{X}})}{\left\|V^{-1/2}X^{(i)}\right\|^2_{\mathcal{X}}} \left(V^{-1/2}X^{(i)}\right)\otimes\left(V^{-1/2}X^{(i)}\right)\\
&& + 2   \left(  \frac{\xi(\left\|V^{-1/2}X^{(i)}\right\|_{\mathcal{X}}) - b_0}{q\beta_1} -    \frac{\psi(\left\|V^{-1/2}X^{(i)}\right\|_{\mathcal{X}})\left\| V^{-1/2}X^{(i)}\right\|_{\mathcal{X}}}{\beta_3}\right) \mathbb{I}  \bigg\} + o_P\left(1\right) \bigg)V^{1/2} \\
&=& -   \frac{1}{  \sqrt{n}}  \sum_{i=1}^{n}\bigg\{  \frac{\theta(\left\|V^{-1/2}X^{(i)}\right\|_{\mathcal{X}})}{\left\|V^{-1/2}X^{(i)}\right\|^2_{\mathcal{X}}} X^{(i)}\otimes X^{(i)}\\
&& + 2   \left(  \frac{\xi(\left\|V^{-1/2}X^{(i)}\right\|_{\mathcal{X}}) - b_0}{q\beta_1} -    \frac{\psi(\left\|V^{-1/2}X^{(i)}\right\|_{\mathcal{X}})\left\| V^{-1/2}X^{(i)}\right\|_{\mathcal{X}}}{\beta_3}\right) V  \bigg\} + o_P\left(1\right)\\
&=& - \widehat{W}_n + o_P\left(1\right)
\end{eqnarray*}
\noindent
where $\widehat{W}_n = n^{-1/2}\sum_{i=1}^{n} \mathcal Z_i ,$  with
\begin{eqnarray*}
\mathcal Z_i &=&   \frac{\theta(\left\|V^{-1/2}X^{(i)}\right\|_{\mathcal{X}})}{\left\|V^{-1/2}X^{(i)}\right\|^2_{\mathcal{X}}} X^{(i)}\otimes X^{(i)} +  \phi\left(\left\|V^{-1/2}X^{(i)}\right\|_{\mathcal{X}}\right)\,V 
\end{eqnarray*}
\noindent
 Slustky's theorem  permits to
conclude that $\widetilde{H}_{n}$  has the same limiting distribution than $-\widehat{W}_n ,$ which can be obtained by using central limit theorem. For doing that$,$ we first have to check that $\mathbb{E}(\mathcal Z_i) = 0$ and $\mathbb{E}\left(\Vert \mathcal Z_i\Vert^2\right)<+\infty$, where $\Vert\cdot\Vert$ is the operator norm induced by te inner product $<A,B>=\textrm{tr}(AB^\ast)$. For proving this  last property we consider the inequality
\begin{eqnarray*}
\mathbb{E}\left(\Vert \mathcal Z_i\Vert^2\right)&\leq&2\mathbb{E}\bigg[\bigg(\frac{\theta(\left\|V^{-1/2}X^{(i)}\right\|_{\mathcal{X}})}{\left\|V^{-1/2}X^{(i)}\right\|^2_{\mathcal{X}}}\bigg)^2\Vert  X^{(i)}\otimes X^{(i)}\Vert^2\bigg]\\
& &+2\Vert V\Vert^2\,\mathbb{E}\bigg( \phi\left(\left\|V^{-1/2}X^{(i)}\right\|_{\mathcal{X}}\right)^2\bigg)\\
&=&2\mathbb{E}\bigg[\bigg(\frac{\theta(\left\|V^{-1/2}X^{(i)}\right\|_{\mathcal{X}})}{\left\|V^{-1/2}X^{(i)}\right\|^2_{\mathcal{X}}}\bigg)^2\Vert  X^{(i)}\Vert^4_\mathcal{X}\bigg]\\
& &+2\Vert V\Vert^2\,\mathbb{E}\bigg( \phi\left(\left\|V^{-1/2}X^{(i)}\right\|_{\mathcal{X}}\right)^2\bigg).
\end{eqnarray*}
Assumption $(\mathscr{A}_4)$ implies that there exists $C>0$ such that $\sup_{t\in\mathbb{R}_+}\left(\theta(t)/t^2\right)\leq C$. On the other hand, from the proof of Corollary 2  in  Lopuha$\ddot{\textrm{a}}$ (1997) (see  p. 236) the functions  $t \mapsto \xi(t)$ and $t \mapsto \psi(t)t$   are \ bounded; then, $\phi$ is also bounded and 
$\mathbb{E}\left( \phi\left(\left\|V^{-1/2}X^{(i)}\right\|_{\mathcal{X}}\right)^2\right)< +\infty$ . Hence
\begin{eqnarray*}
\mathbb{E}\left(\Vert \mathcal Z_i\Vert^2\right)&\leq&2C^2\mathbb{E}\bigg(\Vert  X^{(i)}\Vert^4_\mathcal{X}\bigg)+2\Vert V\Vert^2\,\mathbb{E}\bigg( \phi\left(\left\|V^{-1/2}X^{(i)}\right\|_{\mathcal{X}}\right)^2\bigg)
\end{eqnarray*}
and since $\mathbb{E}\bigg(\Vert  X^{(i)}\Vert^4_\mathcal{X}\bigg)<+\infty$, we deduce that   $\mathbb{E}\left(\Vert \mathcal Z_i\Vert^2\right)<+\infty$. Putting
$Y^{(i)}=V^{-1/2}X^{(i)},$ we have
\begin{eqnarray}\label{ezi}
\mathbb{E}\left(\mathcal Z_i\right) &=& V^{1/2}\mathbb{E}\left(  \frac{\theta\left(\left\|Y^{(i)}\right\|_{\mathcal{X}}\right)}{\left\|Y^{(i)}\right\|^2_{\mathcal{X}}} Y^{(i)}\otimes Y^{(i)} + \phi\left(\left\|V^{-1/2}X^{(i)}\right\|_{\mathcal{X}}\right)\, \mathbb{I} \right)V^{1/2} \ \ \ 
\end{eqnarray}
\noindent
and since $Y^{(i)}$ has a spherical distribution$,$ from equation (4) in Lopuha$\ddot{\textrm{a}}$ (1997) (see  p. 222)   we obtain  
$
\mathbb{E}\left(\xi\left(\left\|Y^{(i)}\right\|_{\mathcal{X}}\right) - b_0\right) = 0$.
\noindent
Further$,$ we have
\begin{eqnarray*}
\mathbb{E}\left(\theta\left(\left\|Y^{(i)}\right\|_{\mathcal{X}}\right)\right) &=& 2q\beta^{-1}_3\mathbb{E}\left(\psi\left(\left\|Y^{(i)}\right\|_{\mathcal{X}}\right)\left\|Y^{(i)}\right\|_{\mathcal{X}}\right).
\end{eqnarray*}
\noindent 
Therefore (\ref{ezi})  becomes:
\begin{eqnarray}\label{ezi2}
\mathbb{E}\left(\mathcal Z_i\right) &=& V^{1/2}\left(\mathbb{E}\left(  \frac{2q\beta^{-1}_3\psi\left(\left\|Y^{(i)}\right\|_{\mathcal{X}}\right)}{\left\|Y^{(i)}\right\|_{\mathcal{X}}} Y^{(i)}\otimes Y^{(i)}\right)  -    \frac{1}{q}\mathbb{E}\left(\theta\left(\left\|Y^{(i)}\right\|_{\mathcal{X}}\right)\right) \mathbb{I} \right)V^{1/2}. 
\end{eqnarray}
\noindent
From Lemma 1 in Lopuha$\ddot{\textrm{a}}$ (1997) (see  p. 221) we have:
\begin{eqnarray*}
\mathbb{E}\left(  \frac{2q\beta^{-1}_3\psi\left(\left\|Y^{(i)}\right\|_{\mathcal{X}}\right)}{\left\|Y^{(i)}\right\|_{\mathcal{X}}} Y^{(i)}\otimes Y^{(i)}\right) &=&   \frac{1}{q}\mathbb{E}\left(  \frac{2q\beta^{-1}_3\psi\left(\left\|Y^{(i)}\right\|_{\mathcal{X}}\right)}{\left\|Y^{(i)}\right\|_{\mathcal{X}}} \left\|Y^{(i)}\right\|^2_{\mathcal{X}}\right)\mathbb{I}\\
&=&   \frac{1}{q}\mathbb{E}\left(2q\beta^{-1}_3\psi\left(\left\|Y^{(i)}\right\|_{\mathcal{X}}\right) \left\|Y^{(i)}\right\|_{\mathcal{X}}\right)\mathbb{I}\\
&=&  \frac{1}{q}\mathbb{E}\left(\theta\left(\left\|Y^{(i)}\right\|_{\mathcal{X}}\right)\right)\mathbb{I}.
\end{eqnarray*}
\noindent
Then, (\ref{ezi2})  implies $\mathbb{E}\left(\mathcal Z_i\right)= 0$. Now$,$ using the central limit theorem we conclude that $-\widehat{W}_n$ 
converges in distribution$,$ as $n \mapsto + \infty,$ to a normal distribution
$N(0, \Lambda)$ in $\mathcal L(\mathcal  X),$ \ where 
$\Lambda$ is the covariance operator of
\begin{eqnarray*}
\mathcal Z &=& -  \frac{\theta(\left\|V^{-1/2}X\right\|_{\mathcal{X}})}{\left\|V^{-1/2}X\right\|^2_{\mathcal{X}}} X\otimes X\\
&& - 2   \left(  \frac{\xi(\left\|V^{-1/2}X\right\|_{\mathcal{X}}) - b_0}{q\beta_1} -    \frac{\psi(\left\|V^{-1/2}X\right\|_{\mathcal{X}})\left\| V^{-1/2}X\right\|_{\mathcal{X}}}{\beta_3}\right) V \\
&=& -2q\beta^{-1}_3\psi\left(\left\|V^{-1/2}X\right\|_{\mathcal{X}}\right)\left\|V^{-1/2}X\right\|^{-1}X\otimes X\\
&& - 2   \left(  \frac{\xi(\left\|V^{-1/2}X\right\|_{\mathcal{X}}) - b_0}{q\beta_1} -    \frac{\psi(\left\|V^{-1/2}X\right\|_{\mathcal{X}})\left\| V^{-1/2}X\right\|_{\mathcal{X}}}{\beta_3}\right) V.
\end{eqnarray*}
\noindent

\subsection{Proof of the theorem}

\noindent
Arguing as in the proof of Theorem 3.2 in Nkiet (2017)  (see p. 203), we have the equality $  \sqrt{n}\left(\widetilde{T}_n - T\right)= \widehat{\varphi}_n\left(\widetilde{H}_{n} \right)$,  where   $\widehat{\varphi}_n$ is the random operator from \ $\mathcal L(\mathcal{X})$   to itself 
defined by:
\begin{eqnarray*}
 \widehat{\varphi}_{n}(A)&=& -f(\widetilde{V}_{n})^{-1}f(A)\left(f(\widetilde{V}_{n})^{-1/2} + \mathbb{I}\right)^{-1}g(\widetilde{V}_{n})f(\widetilde{V}_{n})^{-1/2} + g(A)f(\widetilde{V}_{n})^{-1/2}\\
& &-g(V)(f(\widetilde{V}_{n})^{-1/2}f(A)\left(f(\widetilde{V}_{n})^{-1/2} + \mathbb{I}\right)^{-1}.
\end{eqnarray*}
Considering the linear map $\varphi$  from $\mathcal{X}$ to itself defined as $\varphi(A)=-  \frac{1}{2}f(A)g(V) + g(A) -  \frac{1}{2}g(V)f(A) $
and denoting by \ $\left\|\cdot\right\|_{\infty}$ \  and  \ $\left\|\cdot\right\|_{\infty \infty}$ the  norm  of  \ $\mathcal L(\mathcal{X})$ \  and \ $\mathcal L(\mathcal L(\mathcal{X})),$ respectively defined  by
$\left\|A\right\|_{\infty}=\text{sup}_{x \in \mathcal{X}-\{0\}}\left\|Ax\right\|_\mathcal{X}/\left\|x\right\|_\mathcal{X}$ and $\left\|Q\right\|_{\infty \infty}=\text{sup}_{B \in \mathcal L(\mathcal{X})-\{0\}}\left\|Q(B)\right\|_{\infty}/\left\|B\right\|_{\infty} ,$ we have :
\begin{eqnarray}\label{in1}
\left\|\widehat{\varphi}_n(\widetilde{H}_{n} )-\varphi(\widetilde{H}_{n} )\right\|_{\infty}&\leq&\left\|\widehat{\varphi}_n-\varphi\right\|_{\infty\infty}\left\|\widetilde{H}_{n} \right\|_{\infty} 
\end{eqnarray}
\noindent
and 
\begin{eqnarray}\label{in2}
\left\|\widehat{\varphi}_n-\varphi\right\|_{\infty\infty}&\leq& \Bigg(\left\|f\right\|_{\infty \infty}\left\|\left(f(\widetilde{V}_{ n})^{-1/2} +  \mathbb{I}\right)^{-1}g(\widetilde{V}_{n})f(\widetilde{V}_{n})^{-1/2}\right\|_{\infty}\nonumber\\
&&+ \left\|f\right\|_{\infty \infty}\left\|g(V)\right\|_{\infty} + \left\|g\right\|_{\infty \infty}\nonumber\\
&&+ \left\|f\right\|_{\infty \infty}\left\|g(V)\right\|_{\infty}\left\|\left(f(\widetilde{V}_{n})^{-1/2} +  \mathbb{I}\right)^{-1}\right\|_{\infty} \Bigg)\left\|f(\widetilde{V}_{n})^{-1/2} -  \mathbb{I}\right\|_{\infty}\nonumber\\
&&+ \Bigg(\left\|f\right\|_{\infty \infty}\left\|g(\widetilde{V}_n)f(\widetilde{V}_{n})^{-1/2}\right\|_{\infty}\nonumber\\
&&+ \left\|f\right\|_{\infty \infty}\left\|g(V)\right\|_{\infty}\Bigg)\left\|\left(f(\widetilde{V}_{n})^{-1/2} +  \mathbb{I}\right)^{-1} -   \frac{1}{2}\mathbb{I} \right\|_{\infty} \nonumber\\
&&+   \frac{1}{2}\left\|f\right\|_{\infty \infty}\left\|g\right\|_{\infty \infty}\left\|f(\widetilde{V}_{n})^{-1/2}\right\|_{\infty}\left\|\widetilde{V}_{n} - V \right\|_{\infty} \ \ \ 
\end{eqnarray}
Lemma 1 implies that $\widetilde{V}_n$ converges in probability to 
$V,$ as $n \rightarrow +\infty$.
Then$,$ \ using the \ continuity\\
\noindent
of maps $f$, $g$, $ A\mapsto A^{-1}$ and $A\mapsto A^{-1/2}$ we
deduce that $f(\widetilde{V}_n)$ (resp. $f(\widetilde{V}_n)^{-1}$; resp. $f(\widetilde{V}_n)^{-1/2}$; resp. $g(\widetilde{V}_n)$  converges
in probability$,$ as $n \rightarrow +\infty,$ to \ $\mathbb{I}$ (resp. \ $\mathbb{I};$ resp. \ $\mathbb{I};$ resp. \ $g(V))$. \ 
Consequently 
from (\ref{in1}) and (\ref{in2}) we deduce that $\widehat{\varphi}_n(\widetilde{H}_{n} )-\varphi(\widetilde{H}_{n} )$ converges
in probability to $0$ as $n\rightarrow+\infty$.
Slutsky's theorem allows to conclude that
$\widehat{\varphi}_n(\widetilde{H}_{n} )$ and $\varphi(\widetilde{H}_{n} )$ both converge to the same distribution,  that is the distribution of $\varphi_s\left(M_{s}\right)$. Since $\varphi_s$ is linear this distribution is the normal
distribution with mean equal to $0$ and covariance operator equal to that of
the random variable:
\begin{eqnarray*}
Z_s &=& \varphi\left(\mathcal Z\right) =  \frac{-2q\beta^{-1}_3\psi\left(\left\|V^{-1/2}X\right\|_{\mathcal{X}}\right)}{\left\|V^{-1/2}X\right\|} \varphi\left(X\otimes X\right) - \phi\left(\left\|V^{-1/2}X\right\|_{\mathcal{X}}\right)\,\varphi\left(V\right).
\end{eqnarray*}
\noindent
Besides
\begin{eqnarray*}
\varphi\left(X\otimes X\right)&=&  \sum_{k=1}^K \sum_{\underset{\ell\neq k}{\ell=1}}^K \bigg\{- \frac{1}{2}\bigg(\tau_k^*\left(X_k\otimes X_k\right)  V_{k\ell}\tau_\ell \ + \tau_\ell^*V_{\ell k}\left(X_k\otimes X_k\right)\tau_k\bigg)\\
& & + \tau_k^* \left(X_\ell\otimes X_k\right)\tau_\ell\bigg\} , 
\end{eqnarray*}
and from $f(V) = \mathbb{I},$  it follows:
$
\varphi(V)=g(V) - g(V_s) = g(V) - g(V)= 0$. 
Thus
\begin{eqnarray*}
Z_s &=&  \frac{-2q\beta^{-1}_3\psi\left(\left\|V^{-1/2}X\right\|_{\mathcal{X}}\right)}{\left\|V^{-1/2}X\right\|_{\mathcal{X}}}  \sum_{k=1}^K \sum_{\underset{\ell\neq k}{\ell=1}}^K \bigg\{- \frac{1}{2}\bigg(\tau_k^*\left(X_k\otimes X_k\right)  V_{k\ell}\tau_\ell + \tau_\ell^*V_{\ell k}\left(X_k\otimes X_k\right)\tau_k\bigg)\\
&& + \tau_k^* \left(X_\ell\otimes X_k\right)\tau_\ell\bigg\} .
\end{eqnarray*}
\section{Proof of Theorem 4.2}
\noindent
Under $\mathscr{H}_0$ we have $T = 0$ and,  therefore,  $\sqrt{n}\widetilde{T}_{n} = \sqrt{n}\left(\widetilde{T}_{n} - T\right)$. Consequently,  from Theorem 4.1  we deduce that $\sqrt{n}\widetilde{T}_{n}$ converges in distribution,  as $n \rightarrow +\infty$,  to a random variable $U$ which has a normal distribution in $\mathcal L(\mathcal X)$ with mean 0 and covariance operator equal to that of $Z_s$. Since the map \ $A \mapsto\sum_{k=2}^{K}\sum_{\ell=1}^{k-1}tr\left(\pi_{k\ell}\left(A\right)\pi_{k\ell}\left(A\right)^*\right)$ is continuous$,$ we deduce that $n\widetilde{S}_{n}$ converges in distribution,  as  $n \rightarrow +\infty$,
to
\begin{eqnarray*}
\mathcal Q &=&  \sum_{k=2}^{K}\sum_{\ell=1}^{k-1}\textrm{tr}\left(\pi_{k\ell}\left( U\right)\pi_{k\ell}\left( U \right)^*\right).
\end{eqnarray*} 
\\
\noindent
On the  other hand,    Theorem 4.1 in Nkiet (2017)  shows that  $\mathcal Q= \mathbb{W}^T\mathbb{W}$  where  $\mathbb{W}$ is a random variable having a centered normal distribution in $\mathbb R^d$  with covariance matrix $\Theta$  defined in Nkiet (2017)   with           
\begin{eqnarray*}
\gamma_{ijpt}^{k\ell,ru} &=& < \mathbb{E}\bigg(\pi_{k\ell}\left( U \right)\widetilde{\otimes}\pi_{ru}\left( U \right)\bigg)\left(e_j^{(\ell)}\otimes e_i^{(k)}\right) , e_t^{(u)}\otimes e_p^{(r)} >\\
&=& < \mathbb{E}\bigg(\pi_{k\ell}\left( Z_s \right)\widetilde{\otimes}\pi_{ru}\left( Z_s \right)\bigg)\left(e_j^{(\ell)}\otimes e_i^{(k)}\right) , e_t^{(u)}\otimes e_p^{(r)} >
\end{eqnarray*}
where   $\widetilde{\otimes}$  denotes the tensor product related to the inner product of operators
$<  A , B > = tr (AB^*)$ and \ $\left\{e^{(k)}_i\right\}_{1\leq i \leq p_k}$  is an orthonormal basis of $\mathcal X_k$. Since
under $\mathscr{ H}_0$ we have $V = \mathbb{I}$,  $Z_s$ becomes  
\begin{eqnarray*}
Z_s &=& \displaystyle\frac{-2q\beta^{-1}_3\psi\left(\left\|X\right\|_{\mathcal X}\right)}{\left\|X\right\|_{\mathcal X}}\displaystyle \sum_{k=1}^K \sum_{\underset{\ell\neq k}{\ell=1}}^K \tau_k^* \left(X_\ell\otimes X_k\right)\tau_\ell.
\end{eqnarray*}
\noindent
Thus
\begin{eqnarray*}
\pi_{k\ell}\left( Z_s \right) &=& \frac{-2q\beta^{-1}_3\psi\left(\left\|X\right\|_{\mathcal X}\right)}{\left\|X\right\|_{\mathcal X}}\left(X_\ell\otimes X_k\right)
\end{eqnarray*}
and
\begin{eqnarray*}
\gamma_{ijpt}^{k\ell,ru} &=& -2q\beta^{-1}_3\mathbb{E}\bigg(\displaystyle\frac{\psi\left(\left\|X\right\|_{\mathcal X}\right)}{\left\|X\right\|_{\mathcal X}}< \left(\left( X_\ell\otimes X_k\right)\widetilde{\otimes}\left( X_u\otimes X_r\right)\right)\left(e_j^{(\ell)}\otimes e_i^{(k)}\right) , e_q^{(u)}\otimes e_p^{(r)}>\bigg)\\
&=& -2q\beta^{-1}_3\mathbb{E}\bigg(\displaystyle\frac{\psi\left(\left\|X\right\|_{\mathcal X}\right)}{\left\|X\right\|_{\mathcal X}}< X_\ell\otimes X_k , e_j^{(\ell)}\otimes e_i^{(k)}> < X_u\otimes X_r , e_t^{(u)}\otimes e_p^{(r)}>\bigg)\\
&=&  -2q\beta^{-1}_3\mathbb{E}\bigg(\displaystyle\frac{\psi\left(\left\|X\right\|_{\mathcal X}\right)}{\left\|X\right\|_{\mathcal X}}< X_k , e_i^{(k)}>_k < X_r , e_p^{(r)}>_r < X_\ell , e_j^{(\ell)}>_\ell < X_u , e_t^{(u)}>_u \bigg)\\
&=&  -2q\beta^{-1}_3\mathbb{E}\bigg(z\left(\left\|X\right\|^2_{\mathcal X}\right)< X_k , e_i^{(k)}>_k < X_r , e_p^{(r)}>_r < X_\ell , e_j^{(\ell)}>_\ell < X_u , e_t^{(u)}>_u \bigg),
\end{eqnarray*}
\noindent
where \ $z: t \mapsto \displaystyle\frac{\psi\left(\displaystyle\sqrt{t}\right)}{\displaystyle\sqrt{t}}$. 
Therefore$,$ if $(k, \ell) = (r, u)$ and $(i, j) = (p, t)$ \ with \ $\ell \neq k$ and $u \neq r,$ then
\begin{eqnarray*}
\gamma_{ijpt}^{k\ell,ru} &=& -2q\beta^{-1}_3\mathbb{E}\bigg(z\left(\left\|X\right\|^2_{\mathcal X}\right)< X_k , e_i^{(k)}>_k^2 < X_\ell , e_i^{(\ell)}>_l^2\bigg)
\end{eqnarray*}
\noindent
and  from Lemma 1 in  Lopuha$\ddot{\textrm{a}}$ (1997) we deduce that
\begin{eqnarray*}
\gamma_{ijpt}^{k\ell,ru} &=& \displaystyle\frac{-2q\beta^{-1}_3}{q(q + 1)}\mathbb{E}\bigg(z\left(\left\|X\right\|^2_{\mathcal X}\right)\left\|X\right\|_{\mathcal X}^4\bigg)=\frac{-2\beta^{-1}_3}{(q + 1)}\mathbb{E}\bigg(\psi\left(\left\|X\right\|_{\mathcal X}\right)\left\|X\right\|_{\mathcal X}^3\bigg).
\end{eqnarray*}
\noindent
Otherwise,  if one of the conditions $(k,\ell) = (r, u)$ and $(i, j) = (p, t)$  with \ $\ell \neq k$ and $u \neq r$
does not hold then  $\gamma_{ijpt}^{k\ell,ru} = 0$ . We deduce that
\begin{eqnarray*}
\Theta &=& \displaystyle\frac{-2\beta^{-1}_3}{(q + 1)}\mathbb{E}\bigg(\psi\left(\left\|X\right\|_{\mathcal X}\right)\left\|X\right\|_{\mathcal X}^3\bigg)I_d
\end{eqnarray*}
\noindent
where $I_d$ is the $d\times d$ identity matrix. Thus$,$ $\mathcal Q= \displaystyle\frac{-2\beta^{-1}_3}{(q + 1)}\mathbb{E}\bigg(\psi\left(\left\|X\right\|_{\mathcal X}\right)\left\|X\right\|_{\mathcal X}^3\bigg)\mathcal Q'$ \ where $\mathcal Q'$ is a random variable with distribution equal to $\chi_d^2.$


\begin{thebibliography}{00}






\bibitem{crouxdehon}C. Croux, C.  Dehon, Analyse canonique bas\'ee sur des estimateurs robustes de la matrice des covariances, Rev. Stat. Appl.  50 (2002), 5-26.


\bibitem{crouxfilm}C. Croux,   P. Filmozer,  K. Joosens$,$ Classification efficiencies for robust linear discriminant analysis,  Statistica Sinica 18  (2008), 581-599.



\bibitem{crouxhaes}C. Croux, G.  Haesbroeck G(2000),  Principal component analysis based on robust estimators of the covariance or correlation matrix : Influence function and efficiencies,  Biometrika  87 (2000), 603-618.



\bibitem{davies}P.L. Davies, Asymptotic behaviour of S-estimates of multivariate location and parameters and dispersion matrices, Ann. Statist.  15  (1987),  1269-1292.


\bibitem{gardner}S. Gardner, J.C.  Gower, N.J.  le Roux,  A synthesis of canonical variate analysis, generalized canonical correlation and Procrustes analysis,  Comput. Statist. Data Anal.   50 (2006), 107-134. 

\bibitem{gifi}A. Gifi,  Nonlinear multivariate analysis. Chichesteer : Wiley, 1990.

\bibitem{hampel}F.R. Hampel, E.M.  Ronchetti, P.J.  Rousseeuw,  W.A. Stahel,  Robust Statistics : The Approach based on influence functions. New York :  Wiley, 1986.

\bibitem{hwang}H. Hwang,  K. Jung, Y. Takane, T.S. Woodward,  Functional multiple-set canonical correlation analysis, 
Psychometrika 77 (2012), 753-775.

\bibitem{lopuhaa2}H.P. Lopuha$\ddot{\textrm{a}}$,  On the relation between S-estimator and M-estimator of multivariate location and covariance,
Ann. Statist. 17  (1989), 1662-1683 .

\bibitem{lopuhaa}H.P. Lopuha$\ddot{\textrm{a}}$, Asymptotic expansion of S-estimators of location and covariance,
Statist.  Neerlandica. Volume 51 (1997), 220-237 .



\bibitem{nkiet}G.M. Nkiet, Asymptotics for multiple-set canonical analysis of Euclidean random variables, Math. Methods Statist. 26 (2017), 196-211.



\bibitem{takane}Y.  Takane, H. Hwang, H.  Abdi,  Regularized multiple-set canonical correlation analysis,  Psychometrika   73 (2008), 753-775.


\bibitem{taskinen}S. Taskinen, I.  Koch, H.  Oja,  Robustifying principal component analysis with spatial
sign vectors, Statist. Probab. Lett. 82 (2012), 765-774.



\bibitem{tenenhaus}A. Tenenhaus, M. Tenenhaus, Regularized generalized canonical correlation analysis,  Psychometrika  76 (2011), 257-284.
 




\end{thebibliography}
\end{document}